\documentclass[12pt]{article}
 
\usepackage{babel,euscript}
\usepackage{amssymb,a4wide,doublespace,amsmath,epsf}

\newtheorem{theorem}{Theorem}[section]
\newtheorem{lemma}[theorem]{Lemma}

\newcommand{\Aff}{\mathrm{Aff}}
\newcommand{\RR}{\mathrm{RR}}

\newtheorem{corollary}[theorem]{Corollary}
\newtheorem{example}[theorem]{Example}

\newcommand{\bproof}{\noindent{\it Proof: }}
\newcommand{\eproof}{$\square$ \\}

\newcommand{\sr}{{\rm sr}}

\newcommand{\ep}{\varepsilon}

\newcommand{\C}{{\mathbb C}}
\newcommand{\Z}{{\mathbb Z}}
\newcommand{\N}{{\mathbb N}}

\newcommand{\cK}{{\cal K}}

\newcommand{\R}{{\mathbb R}}
\newcommand{\Cs}{{C$^*$-algebra}}

\newcount\theTime
\newcount\theHour
\newcount\theMinute
\newcount\theMinuteTens
\newcount\theScratch
\theTime=\number\time
\theHour=\theTime
\divide\theHour by 60
\theScratch=\theHour
\multiply\theScratch by 60
\theMinute=\theTime
\advance\theMinute by -\theScratch
\theMinuteTens=\theMinute
\divide\theMinuteTens by 10
\theScratch=\theMinuteTens
\multiply\theScratch by 10
\advance\theMinute by -\theScratch

\def\today{{\number\day\space
 \ifcase\month\or
  January\or February\or March\or April\or May\or June\or
  July\or August\or September\or October\or November\or December\fi
 \space\number\year}}

\newcommand\Afr{{\mathfrak A}}

\newcommand\eqdef{{\;\overset{\mbox{\scriptsize def}}{=}}}

\newcommand\nm[1]{\|#1\|}

\newcommand\Oc{{\mathcal{O}}}

\newcommand\Proj{{\mathrm{Proj}}}

\newcommand\pt{{\tilde p}}

\newcommand\Qc{{\cal Q}}

\newcommand\qt{{\tilde q}}

\newcommand\spec{\mathrm{sp}}

\newcommand\starop[2]{{
  \displaystyle\operatornamewithlimits{\raisebox{-0.5ex}[1.5ex][0ex]{\rm*}}_{#1}^{#2}}}

\newcommand\Tr{{\mathrm{Tr}}}

\begin{document}

\title{Projections in free product C$^*$--algebras, II}
 
\author{Kenneth J.\ Dykema and Mikael R\o rdam}
\date{22 August, 1999}
\maketitle

\begin{abstract}
Let $(A,\varphi)$ be the reduced free product of infinitely many
C$^*$--algebras $(A_\iota,\varphi_\iota)$ with respect to faithful states.
Assume that the $A_\iota$ are not too small, in a specific sense.
If $\varphi$ is a trace then the positive cone of $K_0(A)$ is determined
entirely by $K_0(\varphi)$.
If, furthermore, the image of $K_0(\varphi)$ is dense in $\R$, then $A$ has
real rank zero.
On the other hand, if $\varphi$ is not a trace then $A$ is simple and purely
infinite.
\end{abstract}

\section*{Introduction}

Let $I$ be a set having at least two elements and, for every $\iota\in I$, let
$A_\iota$ be a unital C$^*$--algebra with a state, $\varphi_\iota$, whose GNS
representation is faithful.
Their reduced free product,
\begin{equation}
\label{eq-rfp}
(A,\varphi)=\starop{\iota\in I}{}(A_\iota,\varphi_\iota)
\end{equation}
was introduced by Voiculescu~\cite{Voi:symm} and independently (in a more
restricted way) by Avitzour~\cite{Avi}.
Thus $A$ is a unital C$^*$--algebra with canonical, injective, unital
$^*$--homomorphisms, $\pi_\iota\colon A_\iota\to A$, and $\varphi$ is a state on $A$ such that
$\varphi\circ\pi_\iota=\varphi_\iota$ for all $\iota$.
It is the natural construction in Voiculescu's free probability theory
(see~\cite{VDN}), and
Voiculescu's theory has been vital to the study of
these C$^*$--algebras.

In~\cite{DykRor:proj}, for reduced free product C$^*$--algebras $A$ as
in~(\ref{eq-rfp}), when all the $\varphi_\iota$ are faithful, we investigated
projections in $A$ and the related topic of positive elements in $K_0(A)$.
The behaviour we discovered,  under mild conditions
specifying that the $A_\iota$ are not too small, depended broadly on whether
$\varphi$ is a trace, (i.e.\ on whether all the $\varphi_\iota$ are
traces).
If $\varphi$ is a not trace then by~\cite{DykRor:proj}
$A$ is properly infinite.
It remained open whether $A$ must be purely infinite.
(Some special classes of reduced free product C$^*$--algebras have
in~\cite{DykRor:pi} and~\cite{Dyk:piII} been shown to be purely infinite.)
When $\varphi$ is a trace, then it follows from~\cite{DykRor:proj} that for
every element, $x$, of the subgroup, $G$, of $K_0(A)$ generated by
$\bigcup_{\iota\in I}K_0(\pi_\iota)(K_0(A_\iota))$,  if $K_0(\varphi)(x)>0$
then $x\ge0$ and if $0<K_0(\varphi)(x)<1$
then there is a projection $p\in A$ such that
$x=[p]_0$.
By work of E.\ Germain~\cite{Ger:red}, \cite{Ger:full}, \cite{Ger:amalg}, if
each $A_\iota$ is an amenable C$^*$--algebra then $K_0(A)=G$ and $G$ can be
found from the groups $K_0(A_\iota)$ by using exact sequences, (and by taking
inductive limits if $I$ is infinite);
hence under the hypothesis of amenability, we used Germain's results to give a
complete characterization of
the positive cone of $K_0(A)$ and of its elements corresponding to projections
in $A$.

In the present paper we investigate similar questions for reduced free
product
C$^*$--algebras,~(\ref{eq-rfp}), when $I$ is infinite and when, for infinitely
many $\iota\in I$, there is a unitary, $u\in A_\iota$ such that
$\varphi_\iota(u)=0$.
We show that in this case,
if $\varphi$ is not a trace then $A$ is purely infinite and simple.
If $\varphi$ is a trace then, although we do not know in general if the subgroup
$G$ described above exhausts $K_0(A)$, we nonetheless show that for every
$x\in K_0(A)$, if $K_0(\varphi)(x)>0$ then $x\ge0$; furthermore, we show that
if $x\in K_0(A)$ and if $0<K_0(\varphi)(x)<1$ then there is a projection
$p\in A$ such that $x=[p]_0$.
We also show that if the image of $K_0(\varphi)$ is dense in $\R$ then $A$
has real rank zero.

The {\it real rank} of a C$^*$--algebra, $A$, is denoted $\RR(A)$ and was
invented by L.G.\ Brown and G.K.\ Pedersen~\cite{BroPed:realrank}.
Of particular interest is the case $\RR(A)=0$, which is defined, for a unital
C$^*$--algebra $A$, to mean that
the invertible self--adjoint elements are dense in the set of all self--adjoint
elements of $A$.
If $\varphi$ is a faithful state on an infinite dimensional, simple
C$^*$--algebra $A$, then a necessary condition for 
$\RR(A)=0$ is that there be projections, $p$, in $A$ such that $\varphi(p)$ is
arbitrarily small and positive; hence in particular, the image of
$K_0(\varphi)$ must be dense in $\R$.
We show that this condition is sufficient when $A$ is a reduced free
product of infinitely many algebras, as above, and when $\varphi$ is a trace.
(Moreover, when $\varphi$ is not a trace then $A$ is purely infinite,
so by a result of S.\ Zhang~\cite{Zha}, $A$ has real rank zero.)

\section{Comparison between positive elements and projections}
Most of this section is a reformulation of results from
\cite{Ror:UHFII}. The proof of Theorem \ref{thm1} below is almost identical to
the proof of \cite[Theorem 7.2]{Ror:UHFII}, but the statements of these two
theorems are quite different. 

We recall the notion of comparison of positive elements as introduced by J.\
Cuntz in \cite{Cuntz:simple} and \cite{Cuntz:dimension} (see also
\cite{Ror:UHFII}). Let $A$ be a \Cs{}, and let $a,b$ be positive elements in
$A$. Then $a \lesssim b$ will mean that there exists a sequence $\{x_n\}_{n\in
\N}$ in $A$ with
$$ \lim_{n \to \infty} \|a - x_n b x_n^* \| = 0.$$ If $p,q \in A$ are
projections, then the definition above of $p \lesssim q$ agrees with the usual
definition: $p= vv^*$ and $v^*v \le q$ for some partial isometry $v \in
A$.

If $A$ is unital, and if $\varphi$ is a state on $A$, then define $D_{\varphi}
\colon A^+ \to [0,1]$ by
$$D_{\varphi}(a) = \lim_{\ep \to 0^+} \varphi(f_{\ep} (a)),$$ where $f_{\ep}
\colon \R^+ \to [0,1]$ is the continuous function, which is zero on
$[0,\ep/2]$, linear on $[\ep/2, \ep]$, and equal to 1 on $[\ep, \infty)$.  If
$\varphi$ is a trace, then $D_{\varphi}$ is a dimension function (in the sense
of Cuntz, \cite{Cuntz:dimension}). Notice that $D_{\varphi}(p) = \varphi(p)$
for all projections $p \in A$. 

We shall use the following facts:
\begin{eqnarray} \label{eqn1} f_{2\ep}(a) \le f_{\delta}(f_{\ep}(a)) \le
  f_{\ep/2}(a), \qquad f_{\ep/2}(a)f_{\ep}(a) = f_{\ep}(a),
\end{eqnarray}
when $\ep>0$ and $0 < \delta \le 1/2$, and, consequently,
$D_{\varphi}(f_{\ep}(a)) \le \varphi(f_{\ep/2}(a))$. Moreover, if $0 \le a \le
1$, then $\varphi(a) \le D_{\varphi}(a)$. 
Recall from~\cite{Rfl:sr} that the {\it stable rank} of a unital
C$^*$--algebra $A$ is equal to $1$ if and only if the set invertible elements
of $A$ is dense in $A$.

\begin{lemma} \label{lemma1} Let $A$ be a \Cs{} of stable rank one, let
  $B$ be a hereditary subalgebra of $A$, let $a$ be a positive element in $B$,
  and let $q$ be a projection in $B$ such that $a \lesssim q$. Then for each
  $\ep > 0$ there is a
  projection $p \in B$ such that $f_{\ep}(a) \le p \sim q$.
\end{lemma}

\bproof Observe first that the comparisons $a \lesssim p$ and $p \sim q$ are
independent of whether they are relative to $A$, $B$, or $\tilde{B}$, where
$\tilde{B}$ denotes the \Cs{} obtained by adjoining a unit to $B$. If follows
from~\cite{Rfl:sr} and~\cite{Brown:stiso} that if
$A$ has stable rank one, then so do $B$ and $\tilde{B}$.
By \cite[Proposition 2.4]{Ror:UHFII}, there is for each $\ep
>0$ a unitary $u$ in $\tilde{B}$ with $uf_{\ep}(a)u^* \in q\tilde{B}q \; (= 
qBq)$. Put $p = u^*qu$. Then $p$ is as desired.
\eproof

\begin{lemma} \label{lemma2} Let $A$ be a \Cs{}, let $a$ be a positive
  element in $A$, and let $p$ be a projection in $A$. Then the following are
  equivalent:
\begin{itemize}
\item[{\rm (i)}] $p \lesssim a$,
\item[{\rm (ii)}] $p = xax^*$ for some $x \in A$,
\item[{\rm (iii)}] $p$ is equivalent to some projection in the hereditary
  subalgebra of $A$ generated by $a$.
\end{itemize}
\end{lemma}

\bproof (i) $\Rightarrow$ (ii). If $p \lesssim a$, then $\|p - yay^*\| < 1/2$
for some $y \in A$. Hence $pyay^*p$ is invertible (and positive) in $pAp$, and
therefore $p = zpyay^*pz^*$ for some $z \in pAp$. 

(ii) $\Rightarrow$ (iii). Put $u = xa^{1/2}$. Then $uu^* = p$, and hence $u$ is
a partial isometry. Put $q = u^*u = a^{1/2}x^*xa^{1/2}$. Then $q$ is a
projection in the hereditary subalgebra generated by $a$, and $q \sim p$.

(iii) $\Rightarrow$ (i). Assume $q$ is a projection in the hereditary
subalgebra generated by $a$, and that $q \sim p$. Then there exists an $n \in
\N$ such that $\|q-a^{1/n}qa^{1/n}\| < 1/2$. It follows that
$qa^{1/n}qa^{1/n}q$ and, consequently, $qa^{2/n}q$, are invertible in
$qAq$. Therefore $q = rqa^{2/n}qr^*$ for some $r \in qAq$. This shows that $p
\sim q \lesssim a^{2/n} \lesssim a$.
\eproof

\begin{lemma}[{\cite[Proposition 2.2]{Ror:UHFII}}] \label{lemma3} Let $A$ be a
  unital \Cs{}, let $a,b$ be positive elements in $A$, and let $\ep >0$. If
  $\|a-b\| < \ep$, then $f_{\ep}(a) \lesssim b$.
\end{lemma}

\begin{lemma} \label{lemma4} Let $A$ be a unital \Cs{}, and let $\Afr$ be a
  dense unital *-subalgebra of $A$. Suppose that for each positive element $a
  \in \Afr$ and each $\ep >0$ there is a projection $p \in A$ and $0 < \delta <
  \ep$ such that $f_{\ep}(a) \le p \le f_{\delta}(a)$. Then $\RR(A)=0$.
\end{lemma}

\bproof
To show that $\RR(A)=0$ it will suffice (by~\cite{BroPed:realrank}) to show
that all
self-adjoint elements in the dense $^*$--subalgebra $\Afr$ can be approximated by
invertible self-adjoint elements. 

Let $a$ be a self-adjoint element in $\Afr$, and write $a = a_+ - a_-$. For
each $n \in \N$ find $\delta_n >0$ and projections $p_n,q_n$ in $A$ such that
$$f_{1/n}(a_+) \le p_n \le f_{\delta_n}(a_+), \qquad f_{1/n}(a_-) \le q_n \le
f_{\delta_n}(a_-).$$
Then $p_n \perp q_n$, $p_n a_+ p_n \to a_+$, and $q_n a_- q_n \to a_-$. Set 
$$b_n = (p_n a_+ p_n + {\textstyle {\frac{1}{n}}}p_n)-(q_n a_- q_n +
{\textstyle {\frac{1}{n}}}q_n) + {\textstyle {\frac{1}{n}}}(1-p_n-q_n).$$
Then each $b_n$ is invertible and self-adjoint, and $b_n \to a$. \eproof

\noindent
Let $\cal Q$ be a compact convex subset of the state space of a unital \Cs{}
$A$, and let $\Aff ({\cal Q})$ denote the real vector space of all affine
continuous functions ${\cal Q} \to \R$. Equip this space with the strict
ordering, i.e., $f \ge 0$ if $f =0$ or if $f(\varphi) > 0$ for all $\varphi
\in {\cal Q}$, and, in turn, with the topology induced by this ordering. All
self-adjoint elements $a \in A$ induce an element $\hat{a} \in \Aff ({\cal
Q})$ through the formula $\hat{a}(\varphi) = \varphi(a)$.
We will consider the interval $[0,1]$ of $\Aff({\cal Q})$, defined by
\[
[0,1]=\{f\in\Aff({\cal Q})\mid0\le f\le1\}.
\]

\begin{theorem} \label{thm1} Let $A$ be a unital \Cs{}, let $\cal Q$ be a
  compact convex subset of the state space of $A$, let $\Pi$ be a subset of
  of the set of projections in $A$, and let $\Afr$ be a dense *-subalgebra of
  $A$, which is closed under continuous function calculus (on its normal
  elements). 

Assume that each state in $\cal Q$ is faithful on $\Afr$, and that the following
  comparison properties hold for all positive elements $a \in \Afr$ and for all
  projections $p \in \Pi$:
\begin{itemize}
\item[{\rm ($\alpha$)}] if $D_{\varphi}(a) < \varphi(p)$ for all $\varphi \in \cal
  Q$, then $a \lesssim p$,
\item[{\rm ($\beta$)}] if $\varphi(p) < D_{\varphi}(a)$ for all $\varphi \in
  \cal Q$, then $p \lesssim a$, and
\item[{\rm ($\gamma$)}] the subset of $\Aff(\cal Q)$ induced by $\Pi$ is dense
  in the interval $[0,1]$ of $\Aff({\cal Q})$.
\end{itemize}
It follows that 
\begin{itemize}
\item[{\rm (i)}] if $\sr(A)=1$, then $\RR(A)=0$, and
\item[{\rm (ii)}] if all nonzero projections in $\Pi$ are infinite and full,
  then $A$ is simple and purely infinite. 
\end{itemize}
\end{theorem}

\bproof (i). We show that the conditions in Lemma \ref{lemma4} are
satisfied. So let $a \in \Afr$ be a positive element, and let $\ep >0$. We must
find $0 < \delta < \ep$ and a projection $p \in A$ with $f_{\ep}(a) \le p \le
f_{\delta}(a)$. 

If $\spec(a) \cap (\ep/8,\ep/4) = \emptyset$, then $p=f_{\ep/8}(a)$ and
$\delta = \ep/8$ will be as desired.

Assume now that $\spec(a) \cap (\ep/8,\ep/4) \ne \emptyset$. Then  $0 \le
\varphi(f_{\ep/4}(a)) < \varphi(f_{\ep/8}(a)) \le 1$ for all $\varphi \in \cal
Q$ because each such $\varphi$ is assumed to be faithful on $\Afr$. Since $\Pi$
is dense in the interval $[0,1]$ of $\Aff({\cal Q})$ there is $q \in \Pi$ 
with $\varphi(f_{\ep/4}(a)) < \varphi(q) < \varphi(f_{\ep/8}(a))$ for all
$\varphi \in \cal Q$. By (\ref{eqn1}),
$$D_{\varphi}(f_{\ep/2}(a)) \le \varphi(f_{\ep/4}(a)) < \varphi(q) <
\varphi(f_{\ep/8}(a)) \le D_{\varphi}(f_{\ep/8}(a)).$$ By assumptions
($\alpha$) and ($\beta$) this
implies that $f_{\ep/2}(a) \lesssim q \lesssim f_{\ep/8}(a)$.

{}From Lemma \ref{lemma2} there is a projection $r$ in the hereditary
subalgebra, $B$, generated by $f_{\ep/8}(a)$ such that $q \sim r$. By Lemma
\ref{lemma1} there is a projection $p \in B$ such that $f_{1/2}(f_{\ep/2}(a))
\le p$ (and $p \sim r$). By (\ref{eqn1}), this entails that $f_{\ep}(a) \le p
\le f_{\ep/16}(a)$. The claim is therefore proved with $\delta = \ep/16$.

(ii). If each non-zero hereditary subalgebra of $A$ contains a full element,
then $A$ must be simple. If, moreover, each such hereditary subalgebra
contains an infinite projection, then $A$ is purely infinite and simple (c.f.\
Cuntz' definition of purely infinite simple \Cs s in 
\cite{Cuntz:KOn}). It therefore suffices to show that each non-zero hereditary
subalgebra of $A$ contains an infinite full projection.

Let $B$ be a non-zero hereditary subalgebra of $A$, and let $b$ be a positive
element in $B$ with $\|b\| =1$. Find a positive element $a \in \Afr$ with
$\|a-b\| < 1/2$. Since each $\varphi \in \cal Q$ is faithful, since
$f_{1/2}(a) \ne 0$, and since $\Pi$ is dense in the interval $[0,1]$ of
$\Aff({\cal Q})$,
there is $q \in \Pi$ with $\varphi(q) < \varphi(f_{1/2}(a))$ for all $\varphi
\in \cal Q$. This implies that $\varphi(q) < D_{\varphi}(f_{1/2}(a))$ for all
$\varphi \in {\cal Q}$, and by assumption ($\beta$) we get $q \lesssim
f_{1/2}(a)$. Using Lemma \ref{lemma3} we obtain that $q \lesssim b$, and Lemma
\ref{lemma2} finally implies that there is a projection $p$ in the hereditary
subalgebra of $A$ generated by $b$ (which is contained in $B$, so that $p \in
B$) such that $p \sim q$. Since $q$ is infinite and full, so is $p$, and the
proof is complete.
\eproof

\section{Application to reduced free products of C$^*$--algebras}

Throughout this section, we consider a reduced free product of C$^*$--algebras,
\begin{equation}
\label{eq-rfpinf}
(A,\varphi)=\starop{\iota\in I}{}(A_\iota,\varphi_\iota),
\end{equation}
where $I$ is an infinite set, where each $\varphi_\iota$ is a faithful state
and where for infinitely many $\iota\in I$ there is a
unitary $u\in A_\iota$ with $\varphi_\iota(u)=0$.
It follows from~\cite{Dyk:faith} that $\varphi$ is faithful on $A$.

Avitzour's result~\cite{Avi} gives that $A$ is simple if, for example,
$\varphi$ is a trace.
Indeed, by partitioning the set $I$ into two suitable subsets, $A$ can
be viewed as a reduced free product,
\[
(A,\varphi)=(B_1,\psi_1)*(B_2,\psi_2),
\]
such that there are unitaries, $u\in B_1$ and $v,w\in B_2$ satisfying that
$\psi_1(u)=0=\psi_2(v)=\psi_2(w)$ and that $v$ and $w$ are ${}^*$--free;
hence also $\psi_2(v^*w)=0$.
(Avitzour's result also applies in somewhat more general instances.)
In addition, if $\varphi$ is a trace then by~\cite{DHR:sr} the stable rank of
$A$ is equal to $1$.

The $K_0$-group, $K_0(D)$, of a \Cs{} $D$ is equipped with a
\emph{positive cone} and a \emph{scale} defined respectively by
\begin{eqnarray*}
K_0(D)^+ & = & \{[p]_0 \mid p \in \Proj(D \otimes \cK) \},\\ \Sigma(D) &= &
\{[p]_0\mid p \in \Proj(D) \},
\end{eqnarray*}
where $\Proj(D)$ is the set of projection in $D$, and where $[ \, \cdot \, ]_0
\colon \Proj(D \otimes \cK) \to K_0(D)$ is the canonical map from which $K_0$
is
defined. The positive cone gives rise to an ordering on $K_0(D)$ by $x \le
y$
if $y-x \in K_0(D)^+$, and $x < y$ if $y-x \in K_0(D)^+ \setminus \{0\}$. 
Each (positive) trace $\varphi$ on $D$ induces a positive
group-homomorphism
$K_0(\varphi) \colon K_0(D) \to \R$ which satisfies $K_0(\varphi)([p]_0) =
\varphi(p)$ for $p \in \Proj(D)$, and
$K_0(\varphi)([p]_0)=(\varphi\otimes\Tr_n)(p)$ for
$p\in\Proj(D\otimes M_n(\C))$, where $\Tr_n$ is the (unnormalized) trace
on $M_n(\C)$.
The ordered abelian group
$(K_0(D),K_0(D)^+)$ is
called \emph{weakly unperforated} if $nx >0$ for some $n \in \N$ and some
$x
\in K_0(D)$ implies that $x \ge 0$.

\begin{theorem}
\label{thm-biggy}
Let 
$$ (A, \varphi) = \starop{\iota \in I}{}(A_{\iota}, \varphi_{\iota})$$
be the reduced free product \Cs{}, where each $A_{\iota}$ is a unital
\Cs{},
$\varphi_{\iota}$ is a faithful state on $A_{\iota}$, the index set $I$ is
infinite, and infinitely many $A_{\iota}$ contain a unitary in the kernel
of
$\varphi_{\iota}$. 

If $\varphi$ is a trace (which is the case if all $\varphi_{\iota}$ are
traces), then 
\begin{itemize}
\item[{\rm (i)}] whenever $p,q \in A \otimes M_n(\C)$ are projections such that
$(\varphi \otimes \Tr_n)(p) < (\varphi \otimes\Tr_n) (q)$, it follows that
$p \lesssim q$;
\item[\rm{(ii)}] the positive cone and the scale of $K_0(A)$ are given by
\begin{eqnarray*} 
K_0(A)^+ & = & \{0\} \cup \{ x \in K_0(D) \mid 0 < K_0(\varphi)(x)  \}, \\
\Sigma(A) & = & \{0, 1 \} \cup \{ x \in K_0(D) \mid 0 < K_0(\varphi)(x) < 1
\},
\end{eqnarray*}
and, as a consequence, $(K_0(A),K_0(A)^+)$ is weakly unperforated;
\item[\rm{(iii)}] $\RR(A) = 0$ if and only if $K_0(\varphi)(K_0(A))$ is
dense
in $\R$.
\end{itemize}

If $\varphi$ is not a trace (i.e., if at least one $\varphi_{\iota}$ is not a
trace), then $A$ is simple and purely infinite.
\end{theorem}
\bproof
We consider,
for every finite subset $F\subseteq I$, the C$^*$--subalgebra, $\Afr_F$, of $A$
generated by $\bigcup_{\iota\in F}\pi_\iota(A_\iota)$,
and we let $\Afr=\bigcup_{F\ll I}\Afr_F$, where the union is over all
finite subsets of $I$.
Note that $\Afr$ is a dense, unital
${}^*$--subalgebra of $A$ that is closed under the continuous functional
calculus.

Suppose that $\varphi$ is a trace, let $n\in\N$ and let $p,q\in A\otimes
M_n(\C)$ be projections with
$(\varphi\otimes\Tr_n)(p)<(\varphi\otimes\Tr_n)(q)$.
Using the density of $\Afr$ in $A$ and continuous functional calculus, we find
a finite subset $F$ of $I$ and projections $\pt,\qt\in\Afr_F\otimes M_n(\C)$
such that $\nm{\pt-p}<1$ and $\nm{\qt-q}<1$.
This implies $\pt\sim p$ and $\qt\sim q$.
There are $n^2$ distinct elements
$\iota(1),\iota(2),\ldots,\iota(n^2)\in I\backslash F$
with unitaries
$u_k\in A_{\iota(k)}$ such that $\varphi_{\iota(k)}(u_k)=0$.
Let $B$ be the C$^*$--algebra generated by $\{u_1,u_2,\ldots,u_{n^2}\}$.
Note that $B$ and $\Afr_F$ are free.
Then as in the proof of Proposition~3.3 of~\cite{DykRor:proj}, from the
unitaries $u_1,u_2,\ldots,u_{n^2}$ we can construct a
Haar unitary, $v\in B\otimes M_n(\C)$ such that $\{v\}$ and
$\Afr_F\otimes M_n(\C)$ are $^*$--free (with respect to the tracial state
$\varphi\otimes(\frac1n\Tr_n)$).
Now $\qt\sim v^*\qt v$ and the pair $\pt$ and $v^*\qt v$ is free; moreover,
$(\varphi\otimes\Tr_n)(v^*\qt v)=(\varphi\otimes\Tr_n)(\qt)$.
So by
Proposition~1.1 of~\cite{DykRor:proj}, $v^*\qt v$ is equivalent to a
subprojection, $r$, of $\pt$; hence $q\lesssim p$.
We have thus proved~(i).

The inclusions $\subseteq$ in (ii) are easy consequences of the fact that
$\varphi$ is faithful. Assume $x \in K_0(A)$ and that
$K_0(\varphi)(x)>0$. Since $A$ is unital, there are $n \in \N$ and projections
$p,q \in A \otimes M_n(\C)$ such that $x = [p]_0-[q]_0$. Now,
$$ (\varphi\otimes\Tr_n)(p)-(\varphi\otimes\Tr_n)(q)=K_0(\varphi)(x)>0. $$
Hence, by (i), $q$ is equivalent to a subprojection $\tilde{q}$ of $p$. Thus
$x = [p - \tilde{q}]_0 \in K_0(A)^+$.

Assume next that $x \in K_0(A)$ and that $0 < K_0(\varphi)(x) < 1$. Then, by
the argument above, $x = [p]_0$ for some projection $p \in A \otimes
M_n(\C)$. Let $1_A$ denote the unit of $A$, and let $e \in A \otimes M_n(\C)$
be the diagonal projection whose upper left corner is $1_A$ and with all other
entries equal to 0.
Then
$(\varphi\otimes\Tr_n)(p)=K_0(\varphi)(x)<1=(\varphi\otimes\Tr_n)(e)$.
By (i), this implies that $p$ is equivalent to a
subprojection $\tilde{p}$ of $e$. Hence $x = [\tilde{p}]_0$, and it is easily
seen that $[\tilde{p}]_0 \in \Sigma(A)$.

Finally, to see that $(K_0(A),K_0(A)^+)$ is weakly unperforated, assume that
$x \in K_0(A)$ and that $nx >0$ for some $n \in \N$. Then $K_0(\varphi)(x) =
\frac{1}{n}K_0(\varphi)(nx) > 0$. Hence $x >0$.
We have thus shown~(ii).

Let
\[
\Pi=\bigcup_{F\ll I}\Proj(\Afr_F).
\]
For the set $\Qc$ used in Theorem~\ref{thm1}, we take the singleton
$\{\varphi\}$.
We now show that, regardless of whether $\varphi$ is a trace or not,
conditions~($\alpha$) and~($\beta$) of Theorem~\ref{thm1} hold for every
$p\in\Pi$ and every positive element, $a\in\Afr$.
Given a positive element $a\in\Afr$ and given $p\in\Pi$, there is a finite
subset $F$ of $I$ such
that $a,p\in\Afr_F$.
Let $u\in A_\iota$, for some $\iota\in I\backslash F$, be a unitary such that
$\varphi_\iota(u)=0$.
Then $\{a,p\}$ and $\{u\}$ are $^*$--free with respect to $\varphi$.
Now it follows that $u^*pu$ is a
projection with $\varphi(u^*pu)=\varphi(p)$, and that $u^*pu$ and $a$ are free.
Hence by Lemma~5.3 of~\cite{DykRor:proj} it follows that
$a\lesssim u^*pu$ if $D_\varphi(a)<\varphi(p)$ and
$u^*pu\lesssim a$ if $\varphi(p)<D_\varphi(a)$.
But $u^*pu\sim p$, so ($\alpha$) and ($\beta$) hold.

Suppose now that $\varphi$ is a trace and that the image of $K_0(\varphi)$
is dense in $\R$, and let us show that $\RR(A)=0$.
We will show that $\{\varphi(p)\mid p\in\Pi\}$ is dense in $[0,1]$, which will
imply that condition~($\gamma$) of Theorem~\ref{thm1} holds.
Since the image of $K_0(\varphi)$ is dense in $\R$, the intersection of this
image with $[0,1]$ is dense in $[0,1]$.
By~(ii), it follows that $\{\varphi(p)\mid p\in\Proj(A)\}$ is dense in $[0,1]$.
Since $\Afr$ is dense in $A$, and using continuous functional calculus, we find
for every $p\in\Proj(A)$,
a projection, $\pt\in\Afr$ such that $\varphi(\pt)=\varphi(p)$.
But $\pt\in\Pi$.
We have shown that condition~($\gamma$) of Theorem~\ref{thm1} holds, and we
have already shown that conditions~($\alpha$) and~($\beta$) hold.
Now using the fact that $\sr(A)=1$, we get from
Theorem~\ref{thm1}(i) that $\RR(A)=0$.
This implies one direction of~(iii), but the other direction follows from
more general results.
Indeed,
the image of $K_0(\varphi)$ will be dense in $\R$ if
$A$ contains at least one projection and if $A$ has no minimal
projections.
Both of these conditions hold if $\RR(A)=0$, and if $A$ is simple
and infinite dimensional, as in our case.

Now suppose that $\varphi$ is not a trace, and let us show that $A$ is purely
infinite and simple.
Let $F$ be a finite subset of $I$ such that for at least three distinct
$\iota\in F$ there is a unitary $u\in A_\iota$ satisfying $\varphi_\iota(u)=0$,
and such that for some $\iota\in F$, $\varphi_\iota$ is not a trace.
Then by Theorem~4 of~\cite{DykRor:proj}, the unit is
a properly infinite projection in $\Afr_F$.
Let $\Pi'\subseteq\Pi$ be the set of all full, properly infinite projections in
$\Afr$.
We have already shown that conditions~($\alpha$) and~($\beta$) are satisfied
for every $a\in\Afr$ and every $p\in\Pi'$.
Since~$1$ is a properly infinite projection in some $\Afr_F$, using
Lemma~\ref{lem-propinf} below we get that $\{\varphi(p)\mid p\in\Pi'\}$ is
dense in $[0,1]$, so condition~($\gamma$) is satisfied.
Tautologically, each $p\in\Pi'$ is infinite and full.
Hence by Theorem~\ref{thm1}(ii), $A$ is purely infinite and simple.
\eproof

\begin{lemma}
\label{lem-propinf}
Let $A$ be a unital C$^*$--algebra in which $1$ is properly infinite and let
$\varphi$ be a state on $A$.
Then for every $t\in\R$, $0<t\le1$, there is a projection $p\in A$ such that
$p\sim1$ and $\varphi(p)=t$.
\end{lemma}
\bproof
Using that $1$ is properly infinite, we find isometries,
$v_1,v_2,\ldots$ in $A$ whose range projections are mutually orthogonal.
These generate a unital C$^*$--subalgebra of $A$ isomorphic to the Cuntz
algebra $\Oc_\infty$.
By Cuntz's paper~\cite{Cuntz:KOn}, it follows that if
$p,q\in\Proj(\Oc_\infty)\backslash\{0,1\}$ and if $[p]_0=[q]_0$ in
$K_0(\Oc_\infty)$ then $p$ is homotopic to $q$.
(Indeed, it follows that $p\sim q$ and $1-p\sim1-q$, hence $p$ is unitarily
similar to $q$.
But the unitary group of $\Oc_\infty$ is connected.)

Now let $\ep>0$.
For some $n$ we must have $\varphi(v_nv_n^*)<\ep$;
let $p=v_nv_n^*$.
Then $q\eqdef1-v_n(1-p)v_n^*$ is a projection in $\Oc_\infty$ with
$[q]_0=[1]_0$, and $\varphi(q)>1-\ep$.
Thus there is a continuous path $r_t$ in $\Proj(\Oc_\infty)$ such that
$r_0=p$ and $r_1=q$.
We have that $r_t\sim1$ for all $t$ and
$\{\varphi(r_t)\mid t\in[0,1]\}\supseteq(\ep,1-\ep)$.
\eproof

Let us now state a straightforward application of Theorem~\ref{thm-biggy}
to reduced group C$^*$--algebras.
For a group, $G$, taken with the discrete topology, the reduced group
C$^*$--algebra of $G$, denoted $C^*_{\text{\rm red}}(G)$, is the C$^*$--algebra
generated by the left regular representation of $G$.
The canonical tracial state, $\tau_G$, is the vector state for the
characteristic function of the identity element of $G$.
The following corollary was cited in~\cite{DykHar}, where also a
partial converse was included.

\begin{corollary}
\label{cor-pi}
Let $I$ be an infinite set and let
$$ G=\starop{\iota\in I}{}G_\iota $$
be the free product of nontrivial groups, $G_\iota$.
Suppose that $G$ has finite subgroups of arbitrarily large order.
Then
$$ \RR(C^*_{\text{\rm red}}(G))=0. $$
\end{corollary}
\bproof
We have
$$ (C^*_{\text{\rm red}}(G),\tau_G)\cong\starop{\iota\in I}{}
(C^*_{\text{\rm red}}(G_\iota),\tau_{G_\iota}). $$
If $x$ is a nontrivial element of $G_\iota$ then the left translation operator,
$\lambda_x\in C^*_{\text{\rm red}}(G_\iota)$ is a unitary and $\tau_{G_\iota}(\lambda_x)=0$.
In order to apply Theorem~\ref{thm-biggy}, it is thus sufficient to show that
$C^*_{\text{\rm red}}(G)$ contains projections whose traces (under $\tau_G$) are arbitrarily
small and positive.
But this is clear, since for a finite group $H$, $C^*_{\text{\rm red}}(H)$ contains
projections of trace $1/|H|$.
\eproof

It follows easily from the Kurosh Subgroup Theorem for free products of groups,
(see page~178 of~\cite{LS:comb}), that $G$ has finite subgroups of arbitrarily
high order only if for every positive integer $n$ there is $\iota\in I$ such
that $G_\iota$ has a finite subgroup of order greater than $n$.

\begin{example}
\label{ex-Zns}
If
$$ G=\starop{n=2}\infty(\Z/n\Z), $$
then $C^*_{\text{\rm red}}(G)$ has real rank zero.
Moreover, this C$^*$--algebra is simple, has unique tracial state, has stable
rank one and is not approximately divisible.
\end{example}
\bproof
It has real rank zero by the above corollary.
It is simple and has unique tracial state by Avitzour~\cite{Avi}.
It has stable rank one by~\cite{DHR:sr}.
That it is not approximately divisible follows from the argument in Example~4.8
of~\cite{BKR:apdiv}, because $C^*_{\text{red}}(G)$ is weakly dense in the group
von Neumann algebra $L(G)$, which does not have central sequences.
\eproof

\vspace{.5cm}

\noindent{\sc Department of Mathematics and Computer Science, Odense University,
Campus\-vej 55, 5230 Odense M, Denmark} \\

\noindent{\sl E-mail address:} {\tt dykema@imada.ou.dk} \\
\noindent{\sl Internet home page:}
{\tt www.imada.ou.dk/$\,\widetilde{\;}$dykema/} \\

\vspace{.5cm}

\noindent{\sc Department of Mathematics, University of Copenhagen,
  Uni\-versit\-ets\-park\-en 5, 2100 Copenhagen {\O}, Denmark}\\

\noindent{\sl E-mail address:} {\tt rordam@math.ku.dk}  \\
\noindent{\sl Internet home page:}
{\tt www.math.ku.dk/$\,\widetilde{\;}$rordam/} \\


\begin{thebibliography}{99}

\bibitem{Avi} {\sc D.\ Avitzour,}
  Free products of C$^*$--algebras,
  {\em Trans.\ Amer.\ Math.\ Soc.} {\bf 271} (1982), 423-465.

\bibitem{BKR:apdiv} {\sc B.\ Blackadar, A.\ Kumjian and M.\ R\o{}rdam,}
  Approximately central matrix units and the structure of noncommutative tori
  {\em K--theory} {\bf 6} (1992), 267-284.

\bibitem{Brown:stiso} {\sc L.G.\ Brown,}
  Stable isomorphism of hereditary subalgebras of C$^*$--algebras,
  {\em Pacific J.\ Math.} {\bf 71} (1977), 335-348.

\bibitem{BroPed:realrank} {\sc L.G.\ Brown, G.K.\ Pedersen,}
  C*-algebras of real rank zero,
  {\em J.\ Funct.\ Anal.} {\bf 99} (1991), 131-149.

\bibitem{Cuntz:simple} {\sc J.\ Cuntz,}
  The structure of multiplication and addition in simple C*-algebras, 
  {\em Math.\ Scand.} {\bf 40} (1977), 215--233.

\bibitem{Cuntz:dimension} {\sc J.\ Cuntz,}
  Dimension functions on simple C*-algebras,
  {\em Math.\ Ann.} {\bf 233} (1978), 145--153.

\bibitem{Cuntz:KOn} {\sc J.\ Cuntz,}
  K-theory for certain C*-algebras,
  {\em Ann.\ of Math.} {\bf 113} (1981), 181--197.

\bibitem{Dyk:faith} {\sc K.J.\ Dykema,}
  Faithfulness of free product states,
  {\em J.\ Funct.\ Anal.} {\bf 154} (1998), 223--229.

\bibitem{Dyk:piII} {\sc K.J.\ Dykema,}
  Purely infinite simple C$^*$-algebras arising from free product
  constructions, II, 
  {\it preprint} (1997).  

\bibitem{DHR:sr} {\sc K.J.\ Dykema, U.\ Haagerup, M.\ R\o rdam,}
  The stable rank of some free product C$^*$--algebras
  {\em Duke Math.\ J.} {\bf 90} (1997), 95-121,
  {\em correction} {\bf 94} (1998), 213.

\bibitem{DykHar} {\sc K.J.\ Dykema, P.\ de la Harpe,}
  Some groups whose reduced C$^*$--algebras have stable rank one,
  {\em J.\ Math.\ Pures Appl.} (to appear).

\bibitem{DykRor:proj} {\sc K.J.\ Dykema, M.\ R\o rdam,}
  Projections in free product C$^*$--algebras,
  {\em Geom.\ Funct.\ Anal.} {\bf 8} (1998), 1-16.

\bibitem{DykRor:pi} {\sc K.J.\ Dykema, M.\ R\o rdam,}
  Purely infinite simple C$^*$-algebras arising from free product
  constructions,
  {\em Canad.\ J.\ Math}, {\bf 50} (1998), 323-342.

\bibitem{Ger:red} {\sc  E.\ Germain,}
  $KK$--theory of reduced free product C$^*$--algebras
  {\em  Duke Math.\ J.} {\bf 82} (1996), 707-723.

\bibitem{Ger:full} {\sc  E.\ Germain,}
  $KK$--theory of the full free product of unital C$^*$--algebras
  {\em  J.\ reine angew.\ Math.} {\bf 485} (1997), 1-10.

\bibitem{Ger:amalg} {\sc  E.\ Germain,}
  Amalgamated free product of C$^*$--algebras and $KK$--theory
  {\em Fields Inst.\ Commun.} {\bf 12}, D.\ Voiculescu, ed., (1997), 89-103.

\bibitem{LS:comb} {\sc R.C.\ Lyndon, P.E.\ Schupp,}
  {\em Combinatorial Group Theory},
  Springer--Verlag, 1977.

\bibitem{Rfl:sr} {\sc M.A.\ Rieffel,}
  Dimension and stable rank in the $K$-theory of $C^*$-algebras,
  {\em Proc.\ London Math.\ Soc. (3)}, {\bf 46} (1983), 301--333.

\bibitem{Ror:UHFII} {\sc M.\ R{\o}rdam,}
  On the structure of simple C*-algebra tensored with a UHF-algebra, II,
  {\em J.\ Funct.\ Anal.} {\bf 107} (1992), 255--269.

\bibitem{Voi:symm} {\sc D.\ Voiculescu,}
  Symmetries of some reduced free product $C^*$-algebras,
  {\em Operator Algebras and Their Connections with Topology and Ergodic
  Theory},
  Lecture Notes in Mathematics, vol.\ 1132, Springer-Verlag, 1985, 556--588.

\bibitem{VDN} {\sc D.\ Voiculescu, K.J.\ Dykema, A.\ Nica,}
  {\em Free Random Variables},
  CRM Monograph Series vol.\ 1,
  American Mathematical Society, 1992.

\bibitem{Zha} {\sc S.\ Zhang,}
  Certain C$^*$--algebras with real rank zero and their corona and multiplier
  algebras. I.
  {\em Pacific J.\ Math.} {\bf 155} (1992), 169--197.

\end{thebibliography}
\end{document}